\documentclass[draft]{amsart}

\usepackage{pgf,tikz, ytableau}
\usetikzlibrary{arrows}
\usepackage{cite}

\usepackage{amssymb,graphicx,enumerate,color, hyperref}
\usepackage{multirow}

\renewcommand{\l}{\lambda}

\renewcommand{\P}{\mathcal{P}}

\newcommand{\D}{{\mathcal{D}}}

\newcommand{\R}{{\mathcal{R}}}

\newcommand{\beqs}{\begin{equation}}
\newcommand{\eeqs}{\end{equation}}
\numberwithin{equation}{section}
 \theoremstyle{plain}
\newtheorem{theorem}{Theorem}[section]

\newtheorem{corollary}[theorem]{Corollary}

\theoremstyle{remark}

\theoremstyle{empty}

\begin{document}

\makeatletter
\def\imod#1{\allowbreak\mkern10mu({\operator@font mod}\,\,#1)}
\makeatother

\author{Ali Kemal Uncu}
   \address{University of Bath, Faculty of Science, Department of Computer Science, Bath, BA2 7AY, UK}
   \email{aku21@bath.ac.uk}

%\scalebox{.9}
\title{Combinatorial Proofs of Ismail's Identities on Al-Salam--Chirara Polynomials}

\dedicatory{Dedicated to Mourad Ismail, a great inspiration,  on the occasion of his 80th birthday.}

\begin{abstract} We present new proofs of two identities arising in work of Mourad Ismail using partition theoretic generating function interpretations.
\end{abstract}   

 %%%%	05A15  	Exact enumeration problems, generating functions [See also 33Cxx, 33Dxx]
 %	05A17  	Partitions of integers [See also 11P81, 11P82, 11P83]
 %  05A19  	Combinatorial identities, bijective combinatorics
 % 	05A20  	Combinatorial inequalities
 %  11B65  	Binomial coefficients; factorials; $q$-identities [See also 05A10, 05A30]
 %	11B34  	Representation functions
 %	11B75  	Other combinatorial number theory
 %	11P81  	Elementary theory of partitions [See also 05A17]
 %	11P84  	Partition identities; identities of Rogers-Ramanujan type
 %	11P99  	None of the above, but in this section
 %	33D15  	Basic hypergeometric functions in one variable, ${}_r\phi_s$

\date{\today}
   
\maketitle

\section{Introduction}

Ismail has been a big force in the progress of $q$-series research. He is a virtuoso of orthogonal polynomials and special functions. I met ---if I may--- Mourad while I was as a graduate student at the University of Florida. Although on paper we were close by, it was seldom that we visited each other. The same applies to research. On paper we were interested in similar mathematical objects, some results we found always looked closely related, yet $q$-series have too many facets and I could not understand how he was making his discoveries. Mourad's expertise was, and still is, magical to me. I was mastering combinatorial treatment of $q$-series through partition theory. I saw more of Mourad in Austria than I saw him in Florida. He visited Austria multiple times, while I was a postdoctoral researcher at the Research Institute for Symbolic Computation and later at the Johann Radon Institute for Computational and Applied Mathematics. In 2018, I watched a talk by Mourad in Peter Paule's (another master virtuoso) 60th birthday celebrations \cite{MouradSlides}. We exchanged some thoughts about the talk and the contents, especially about combinatorial implications of his results. Mourad pointed me to his (then) manuscript \cite{Mourad2019}. He shared his desire to see if it would be possible to prove his results using combinatorics. I proved some of them and shared these proofs in private conversation. I am grateful that he acknowledged these proofs in \cite{Mourad2019}. It is a badge of honour for an early postdoc to receive such recognition from a mathematics legend. For his 80th birthday, I decided to revisit and rewrite these proofs and dedicate those to him. Hopefully, this will also encourage other partition theory researchers to look into Mourad's and all associated works. It is certain that there is a lot to learn, a lot to prove, and a lot of room to discover and refine on this facet of $q$-series using combinatorial techniques.

In this paper, we present new proofs for the following two theorems using partition generating function interpretations. Let $q$, $a$, $b$, $t$, and $x$ all be formal transcendental variables for the rest of the paper which do not satisfy any relation among themselves.

\begin{theorem}\label{Thm1}
\begin{equation}\label{Ismail1}\frac{(atq;q)_\infty}{(q;q)_\infty}\sum_{k\geq 0} \frac{(-bt)^k}{(q,atq;q)_k}q^{k(k+1)/2} = \frac{(atq,btq;q)_\infty}{(q;q)_\infty}\sum_{k\geq 0} \frac{(-1)^k (abt^2 q)^k q^{k(3k-1)/2}}{(atq,btq,q;q)_k}.\end{equation} 
\end{theorem}

\begin{theorem}\label{Thm2}
	\begin{equation}\label{Ramanujan_F_raw}\sum_{n\geq 0} \frac{q^{n^2} (-x)^n}{(q;q)_n} = (xq^2;q^2)_\infty \sum_{n\geq0} \frac{q^{n^2}(-x)^n}{(q^2,xq^2;q^2)_n}. \end{equation}
\end{theorem}

It should be noted that the above theorems appear in Ramanujan's lost notebooks. Precisely, these are \cite[Entry 2.2.4]{RamaLost2} and \cite[Entry 1.5.1]{RamaLost2}, respectively. A different and much shorter combinatorial proof of Theorem~\ref{Thm2} was given by Berndt--Kim--Yee in \cite[Theorem 5.5]{BerndtCombo}.

I was not aware of \cite{BerndtCombo} when addressing Mourad’s request. The original note, of my proofs to Theorems~\ref{Thm1} and \ref{Thm2}, was written in good length for Mourad's eyes only and assumed no prior partition theory knowledge from him. Nevertheless, independently studying proven identities has once again shown its merits. We exemplify it with the following partition relation and the following $q$-series identities. Different dissections of partitions, that produce new proofs, lead to new bijections and additional combinatorial connections. For example, while proving Theorem~\ref{Thm2}, we will also obtain the following new partition theorem.

\begin{theorem}\label{Thm_Ramanujan_F_combinatorial}
	\[ \sum_{\pi\in\R\R} x^{\nu(\pi)} q^{|\pi|}= \sum_{\substack{\pi\in\D\\ s_e(\pi) > 2 \nu_o(\pi)}} x^{\nu(\pi)} q^{|\pi|},\]
	where $\D$ and $\R\R$ are sets of partitions with gap between parts $\geq 1$ and $\geq 2$, respectively, $\nu(\pi)$ and $\nu_o(\pi)$ are the number of parts and number of odd parts of a partition, respectively, and $s_e(\pi)$ is the smallest even part (or $\infty$ if there are no even parts).
\end{theorem}

To further exemplify the combinatorial connection Theorem~\ref{Thm_Ramanujan_F_combinatorial} yields, we explicitly write the partitions of 15 into exactly 3 parts for both sets in Table~\ref{Table_RR_example}.

\begin{table}[htp]\caption{Partitions of 15 into exactly 3 parts that are counted by the two sides of Theorem~\ref{Thm_Ramanujan_F_combinatorial} }\label{Table_RR_example}
	\[\begin{array}{c|c}
		\pi\in\R\R & \pi\in\D,\ s_e(\pi) > 2 \nu_o(\pi)\\[-1.5ex]\\
		(11,3,1) & (11,3,1)\\
		(10,4,1) & (10,4,1)\\
		(9,5,1) & (9,5,1) \\
		(9,4,2) & (8,6,1)\\
		(8,6,1) & (8,4,3) \\
		(8,5,2) & (7,5,3) \\
		(7,5,3) & (6,5,4) \\
	\end{array}\]
\end{table}

Moreover, understanding the combinatorial structure of the generating functions lead us to refinements of these functions and their identities. The next polynomial identity is one of such identities, which implies Theorem~\ref{Thm1}.

\begin{theorem}\label{Thm4}
	For all non-negative integer $L$, we have
	\begin{align}\nonumber \sum_{n,k=0}^{\lfloor L/2\rfloor} \frac{q^{n^2+2nk+k(k+1)}x^{n+k}}{(q^2;q^2)_{n+k}}{n+k\brack k}_{q^2} & (q^{2(\lfloor (L+1)/2 \rfloor -n-k+1)};q^2)_n (q^{2(\lfloor L/2 \rfloor -n-k+1)};q^2)_k\\\label{thm4_eq} &= \frac{(-1)^L-1}{2}q^{L^2}x^L+ \sum_{n=0}^{L+1} q^{n^2}x^n {L-n+1\brack n}_q, \end{align}
	where  $\lfloor \cdot \rfloor$ is the floor function which yields the greatest integer that is less than or equal to the argument.
\end{theorem}

A corollary of Theorem~\ref{Thm4}, by comparing the exponents of $x$ in \eqref{thm4_eq}, is the following.

\begin{corollary}
	For all non-negative integer $L$ and $m$, we have
	\[ \sum_{k\geq 0} q^k \frac{(q^{2(\lfloor (L+m)/2 \rfloor-m+1)};q^2)_{m-k} (q^{2(\lfloor (L+m-1)/2 \rfloor -m+1)};q^2)_k}{(q^2;q^2)_m} {m\brack k}_{q^2} =  {L\brack m}_q .\]
\end{corollary}

The rest of the paper is organized as follows. The next section has the necessary background information. The following two sections lay out the combinatorial proofs of the Theorems~\ref{Thm1} and \ref{Thm2}, and the observed partition identities. The last section is dedicated to an outlook, a discussion of Theorem~\ref{Thm4} and a similar result that implies Theorem~\ref{Thm1}, and possible future directions.

\section{Some Background}\label{Sec_back}

We will follow classical partition theory definitions of Andrews' encyclopedia \cite{A}. For completeness, we are reciting some important definitions and formulas here.

A partition is a non-increasing finite list of positive integers, such as $\pi= (\lambda_1,\lambda_2,\dots,\lambda_{\nu(\pi)})$, where $\lambda_i$ are the parts of $\pi$, $\nu(\pi)$ denotes the number of parts of $\pi$, and $l(\pi)$ denotes the largest part of $\pi$. The sum of all parts of $\pi$, denoted by $|\pi|$, is called the size of $\pi$. We call a $\pi$ ``\textit{a partition of n}" if the size is equal to $n$. Let $\D$ be the set of partitions into distinct parts and $\P$ be the set of all partitions. For example, $(5,4,1)$ is a partition (into distinct parts) of 10 that has 3 parts. Also the empty list is the unique partition of size 0 (that has 0 parts).

There is a one-to-one correspondence between partitions and Ferrers diagrams (or Young tableaux). For example, the partition $(4,3,3,2)$ is associated with the diagram,
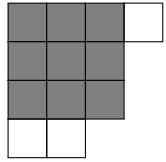
\begin{center}
\ytableausetup{nosmalltableaux, boxsize=5mm}
\begin{ytableau}
	*(gray) {}  & *(gray) {} & *(gray) {} & {}\\
 *(gray)	{} & *(gray) {} & *(gray)  {} \\
*(gray) {}	& *(gray) {} & *(gray) {}\\
{}	&{}
\end{ytableau}
\end{center}
where the largest part of the partition corresponds to the number of boxes on the top row, the second largest part corresponds to the number of boxes on the second row from the top, etc. Finally, there is a unique largest square that can be fit into the Ferrers diagram of a partition where one corner of the square aligns with the top left corner of the diagram. This unique square is called the Durfee square and we will denote an edge length of the Durfee square of $\pi$ as $d(\pi)$. In the above example, the Durfee square is of size $3\times 3$ and shaded gray in the Ferrers diagram.

We define the $q$-Pochhammer symbol \[
(a;q)_n := \prod_{i=0}^{n-1} (1-a q^i)\ \ \text{and}\ \  (a_1,a_2,\dots,a_k;q)_n := \prod_{i=1}^{k} (a_i;q)_n, 
\]
where $n \in \mathbb{Z}_{\geq 0} \cup \{\infty\}$. Also, we define Gaussian coefficients (or $q$-binomial coefficients) as \[{a+b \brack b}_q = \left\{\begin{array}{cc}
	\frac{(q;q)_{a+b}}{(q;q)_a (q;q)_b}, & \text{if } a,b\geq 0,\\
	0& \text{otherwise.}
\end{array}\right.\] The classical partition theoretic interpretation of the Gaussian coefficient \[{a+b\brack b }_q\] is the generating function for the number of partitions where the largest part is $\leq a$ and the number of parts $\leq b$ (equivalently, the Ferrers diagrams fit in an $a \times b $ box).

Let $1\leq j$ and $j\leq N \leq \infty$. The generating function for the number of partitions with parts from the set $\{j, j+1, \dots, N\}$, where the exponent of $q$ keeps track of the size of the partition and the exponent of $x$ keeps track of the number of parts in the partition, can be given in multiple ways:
\begin{equation}
	\label{GF_ordinary_bdd_N}
	\frac{1}{(x q^j;q)_{N-j+1}} = 1+\sum_{i=j}^N \frac{x q^i}{(x q^i;q)_{N-i+1}} = \sum_{i=j}^N \frac{x^i q^{i^2}}{(x q^j;q)_{i-j+1}}{N\brack i}_q + \sum_{i=0}^{j-1} x^i q^{ i j } {N-j+i\brack i}_q.
\end{equation}
These interpretations of the expressions from left-to-right follow from the geometric series expansion and formal multiplications, by dissecting a partition with respect to its smallest part and remaining parts, and dissecting a partition with respect to the Durfee squares, respectively. 
%One can also observe that reflecting a partition (its Ferrers diagram) over the diagonal (from top left to bottom right) of its Durfee square yields another partition. This is called the \textit{conjugate partition}. By looking at the conjugates, onecan also interpret \eqref{GF_ordinary_bdd_N} as the generating function for the number of partitions where the number of parts is in $\in \{j, j+1, \dots, N\}$, where the exponent of $q$ keeps track of the size of the partition and the exponent of $x$ keeps track of the number of parts in the partition,
Similarly, \begin{equation}
	\label{GF_ordinary_bdd_N_parts}
	 \sum_{i=j}^N \frac{x^i q^i}{( q;q)_i}% = \sum_{i=j}^N \frac{x^i q^{i^2}}{(q,x q;q)_i}.
\end{equation}
can be interpreted as the generating function for the number of partitions with at least $j$ parts and at most $N$ parts, where the exponent of $q$ keeps track of the size of the partition and the exponent of $x$ keeps track of the number of parts in the partition. 

One can also write down similar generating functions for partitions into distinct parts. The generating function for the number of partitions into distinct parts from the set $\{j, j+1, \dots, N\}$, where the exponent of $q$ keeps track of the size of the partition and the exponent of $x$ keeps track of the number of parts in the partition, can be represented by $(-x q^j;q)_{N-j+1}$. The generating function for the number of distinct partitions into at least $j$ and at most $N$ parts is
\begin{equation}
	\label{GF_distinct_bdd_N}
	\sum_{i=j}^N \frac{x^i q^{i(i+1)/2}}{(q;q)_i}.
\end{equation}

Here the exponent ${i(i+1)/2}$ is interpreted as the partition $(i,i-1,i-2,\dots,1)$, which is the smallest possible partition into $i$ distinct parts. 

In all of the sum representations of \eqref{GF_ordinary_bdd_N}, \eqref{GF_ordinary_bdd_N_parts}, and \eqref{GF_distinct_bdd_N} the summands represent the partitions with the respective dissection with a fixed $i$, whether it is the number of parts or the side-length of the Durfee square. Summing over all possibilities from $j$ to $N$ yields the generating function of the interpretation. %One can easily see this by picking $j=N=i$ in the respective interpretations.

\section{Proof of Theorem~\ref{Thm1}}

Before we start proving \eqref{Ismail1}, note that the sides of this equation are the right-hand sides of (5.2) and (5.7) in \cite{Mourad2019}. Ismail proves this identity by studying the limiting images of his central object Al-Salam--Chirara polynomials and Jackson's ${}_2\phi_1-{}_2\phi_2$ transformation \cite[III.4, p. 359]{GasperRahman} (which Ismail refers to as the $q$-Pfaff--Kummer transformation).

Here we would like to provide a quick proof using combinatorial reasoning only. After some basic cancellations and $b\mapsto-b$, \eqref{Ismail1} is equivalent to the identity \begin{equation}\label{Ismail1_simplified}\sum_{k\geq 0 }\frac{(bt)^kq^{k(k+1)/2}}{(q,atq;q)_k} = \sum_{m\geq 0}(bt)^mq^{m(m+1)/2}(-btq^{m+1};q)_\infty \frac{(at)^mq^{m^2}}{(q,atq;q)_m}.\end{equation}
Our aim is to show that both sides of the identity are the generating functions for the number of the same set of partition pairs.

Starting from the left-hand side of \eqref{Ismail1_simplified}, for a fixed $k$, we know that \[\frac{(bt)^kq^{k(k+1)/2}}{(q;q)_k}\] is the generating function for the number of partitions into $k$ distinct parts there the exponent of $(bt)$ counts the number of parts (see Section~\ref{Sec_back}). We also know that \[\frac{1}{(atq;q)_k}\] is the generating function for the number of partitions with parts $\leq k$, where the exponent of $(at)$ counts the number of parts. Together, the summand of the left-hand side of \eqref{Ismail1_simplified} for a fixed $k$ is the generating function for the number of partition pairs, where the first partition $\pi_1$ is into exactly $k$ distinct parts and the second partition $\pi_2$ (is not necessarily into distinct parts) satisfies the relation that $l(\pi_2) \leq k$. Summing over all $k$, we see the relation: \begin{equation}
\label{Ismail_LHS}
\sum_{k\geq 0 }\frac{(bt)^kq^{k(k+1)/2}}{(q,atq;q)_k} = \sum_{\substack{(\pi_1,\pi_2)\in \D\times\P\\l(\pi_2)\leq \nu(\pi_1)}} t^{\nu(\pi_1)+\nu(\pi_2)}a^{\nu(\pi_2)}b^{\nu(\pi_1)}q^{|\pi_1|+|\pi_2|}.\end{equation}

Now we move on to the interpretation of the right-hand side of \eqref{Ismail1_simplified}. The expression \[(bt)^m q^{m(m+1)/ 2}(-btq^{m+1};q)_\infty\] is the generating function for the number of partitions $\pi^*_1$ into $\geq m$ distinct parts, where parts $1,2,\dots, m$ all appear in the partition, and the exponent of $(bt)$ keeps track of the number of parts. The fraction \[ \frac{(at)^mq^{m^2}}{(q,atq;q)_m}\] is the generating function for the number of partitions with a $m\times m$ Durfee square, and the exponent of $(at)$ keeps track of the number of parts. Summing over all $m$ one gets 
\begin{equation}
\label{Ismail_RHS}
\sum_{m\geq 0}(bt)^mq^{m(m+1)/ 2}(-btq^{m+1};q)_\infty \frac{(at)^mq^{m^2}}{(q,atq;q)_m} = \sum_{\substack{(\pi^*_1,\pi^*_2)\in\D\times\P\\d(\pi^*_2)\leq r(\pi^*_1)}} t^{\nu(\pi^*_1)+\nu(\pi^*_2)}a^{\nu(\pi^*_2)}b^{\nu(\pi^*_1)}q^{|\pi^*_1|+|\pi^*_2|},\end{equation} where $r(\pi^*_1)$ is the number of consecutive parts in $\pi^*_1$ starting from 1.

Our next task is to show that \eqref{Ismail_LHS} and \eqref{Ismail_RHS} count the same amount. To this end, we would like to modify the left side of \eqref{Ismail_RHS} and show that it can be viewed as a generating function which counts the same partitions counted in \eqref{Ismail_LHS}.

For a fixed $m$, we have \[(bt)^m q^{m(m+1)/ 2}(-btq^{m+1};q)_\infty = \sum_{k=m}^\infty \frac{(bt)^kq^{k(k+1)/2}}{(q;q)_{k-m}}.\] Here the right-hand side summand can be viewed as the generating function for the number of partitions into $k$ distinct parts, where $1,2,\dots, m$ appear as parts, and the exponent of $(bt)$ keeps track of the number of parts. Summation over $k$ on the right give the equivalence of both sides. Then, writing this in the right-hand side of \eqref{Ismail1_simplified}, we get
\begin{align}
\nonumber\sum_{m\geq 0}(bt)^mq^{m+1\choose 2}(-btq^{m+1};q)_\infty \frac{(at)^mq^{m^2}}{(q,atq;q)_m} &= \sum_{m\geq 0}\left(\sum_{k=m}^\infty \frac{(bt)^kq^{k(k+1)/2}}{(q;q)_{k-m}}\right) \frac{(at)^mq^{m^2}}{(q,atq;q)_m},\\
\nonumber&= \sum_{m\geq 0}\sum_{k=m}^\infty \frac{(bt)^kq^{k(k+1)/2}}{(q;q)_{k-m}(q;q)_m} \frac{(at)^mq^{m^2}}{(atq;q)_m},\\
\label{Main_Connection}& = \sum_{m\geq 0}\sum_{k=m}^\infty \frac{(bt)^kq^{k(k+1)/2}}{(q;q)_{k}} {k\brack m }_q \frac{(at)^mq^{m^2}}{(atq;q)_m}.
\end{align} 
%where the Gaussian coefficient \[{k\brack m }_q := \left\{\begin{array}{cc}
%\frac{(q;q)_m}{(q;q)_{k-m}(q;q)_m}, & \text{if }k\geq m\geq0,\\
%0, & \text{otherwise}. 
%\end{array} \right.\] The classical partition theoretic interpretation of the Gaussian coefficient \[{n+m\brack m }_q\] is the generating function for the number of partitions that fit in a $m \times n $ box.

For a fixed $m$, the inner sum  of \eqref{Main_Connection}: \[\sum_{k=m}^\infty \frac{(bt)^kq^{k(k+1)/2}}{(q;q)_{k}} {k\brack m }_q \frac{(at)^mq^{m^2}}{(atq;q)_m}\] can be interpreted in three pieces. First, \[\sum_{k=m}^\infty \frac{(bt)^kq^{k(k+1)/2}}{(q;q)_{k}}\] is the generating function for the number of partitions $\pi_1$ into $k (\geq m)$ distinct parts, where the exponent of $(bt)$ keeps track of the number of parts. Secondly, \[ \frac{(at)^mq^{m^2}}{(atq;q)_m}\] is the generating function for the number of partitions $\hat{\pi}_2$ where the largest part is exactly $m$ and the largest part repeats at least $m$ times, and the exponent of $(at)$ keeps track of the number of parts. Lastly \[{k \brack m}_m\] is the generating function for the number of partitions $\hat{\pi}_3$ that fit into a $m \times (k-m)$ box.

Note that the partitions $\hat{\pi}_2$ and $\hat{\pi}_3$ can be put together to make a new partition $\pi_2$ by placing $\hat{\pi}_3$ on the right of the $m$-many $m$'s. The Durfee square of the resulting partition is the same as the Durfee square of $\hat{\pi}_2$. This operation is bijective. One can revert this action by splitting the parts to the right of the Durfee square for any given partition. This merging operation changes the size of the largest parts.  The partition $\pi_2$ has the largest part $\leq k$ (which is less than the number of parts of $\pi_1$) and its Durfee square is of the size $m\times m$. 
Combining what we have observed, we have \[\sum_{k=m}^\infty \frac{(bt)^kq^{k(k+1)/2}}{(q;q)_{k}} {k\brack m }_q \frac{(at)^mq^{m^2}}{(atq;q)_m} =  \sum_{\substack{(\pi_1,\pi_2)\in \D\times\P\\d(\pi_2)=m\\  l(\pi_2)\leq \nu(\pi_1)}} t^{\nu(\pi_1)+\nu(\pi_2)}a^{\nu(\pi_2)}b^{\nu(\pi_1)}q^{|\pi_1|+|\pi_2|}.\]

Summing over all $m\geq 0$ lifts the restrictions on the Durfee squares' size of the partition $\pi_2$. This shows that the left-hand side of \eqref{Ismail_RHS} is equal to the right-hand side of \eqref{Ismail_LHS}: \begin{equation}\label{Ismail1RHS}\sum_{m\geq 0}(bt)^mq^{m+1\choose 2}(-btq^{m+1};q)_\infty \frac{(at)^mq^{m^2}}{(q,atq;q)_m}=\sum_{\substack{(\pi_1,\pi_2)\in \D\times\P\\l(\pi_2)\leq \nu(\pi_1)}} t^{\nu(\pi_1)+\nu(\pi_2)}a^{\nu(\pi_2)}b^{\nu(\pi_1)}q^{|\pi_1|+|\pi_2|}.\end{equation} Note that the right-hand sides of \eqref{Ismail_LHS} and \eqref{Ismail1RHS} are the equal. This proves \eqref{Ismail1_simplified} and implies Theorem~\ref{Thm1}.

The considerations in this section also proved the following result.

\begin{theorem}\label{Thm3}
	\[\sum_{\substack{(\pi_1,\pi_2)\in \D\times\P\\l(\pi_2)\leq \nu(\pi_1)}} t^{\nu(\pi_1)+\nu(\pi_2)}a^{\nu(\pi_2)}b^{\nu(\pi_1)}q^{|\pi_1|+|\pi_2|} = \sum_{\substack{(\pi^*_1,\pi^*_2)\in\D\times\P\\d(\pi^*_2)\leq r(\pi^*_1)}} t^{\nu(\pi^*_1)+\nu(\pi^*_2)}a^{\nu(\pi^*_2)}b^{\nu(\pi^*_1)}q^{|\pi^*_1|+|\pi^*_2|}.\]
\end{theorem}

The generating functions have a subtle change, but this change is non-trivial. For example, let's give two partition pairs with size total 8. The partition pair $( (3,2),(2,1))$ is counted by the left-hand side sum but not by the right-hand side, and $( (2,1),(3,2))$ is counted by the right and not the left.

\section{Proof of Theorem~\ref{Thm2}}\label{S2}

Now we move onto the task of proving Theorem~\ref{Thm2}, which proves to be slightly harder. Our proof is different from the one in \cite{BerndtCombo}. There the authors interpret both sides of the equation as 2 modular diagrams and argue that both sides have the same number of partitions by plucking out 1 from each odd part and then arguing that the numbers of leftover partitions generated by both sides are the same. We also use $2$-modular diagrams, but in a different way. Here we will define an explicit bijection between the minimal configurations and Ferrers diagrams that fit in a box of certain dimensions, which has $q$-binomial coefficients as their the generating functions. We also prove the partition identity Theorem~\ref{Thm_Ramanujan_F_combinatorial} here, using generating function interpretations.
%, is to give a combinatorial explanation to the identity (5.13): \begin{equation}\label{Ramanujan_F_raw}\sum_{n\geq 0} \frac{q^{n^2} (-x)^n}{(q;q)_n} = (xq^2;q^2)_\infty \sum_{n\geq0} \frac{q^{n^2}(-x)^n}{(q^2,xq^2;q^2)_n}. \end{equation}

The left-hand side of \eqref{Ramanujan_F_raw} is, what Mourad \cite{MouradSlides} refers to as, the Ramanujan function with $x\mapsto -x$. This function is widely seen in combinatorial treatment of partitions. %To our pleasant surprise, this identity can be proven using only combinatorial reasoning as above. 
Before presenting that proof, we would like to give combinatorial interpretations of both sides of this identity. In \eqref{Ramanujan_F_raw} let $x\mapsto -x$ and on the right-hand side carry the outside factor $(-xq^2;q^2)_\infty$ inside and do the simple cancellations that appear in the summand. This yields the equivalent formula \begin{equation}\label{Ramanujan_F}\sum_{n\geq 0} \frac{q^{n^2} x^n}{(q;q)_n} =\sum_{n\geq0} \frac{q^{n^2}x^n}{(q^2;q^2)_n} (-xq^{2n+2};q^2)_\infty . \end{equation} 

The left-hand side of \eqref{Ramanujan_F} is the Ramanujan function and we would like to start our presentation with its classical interpretation. Let $\R\R$ be the set of (Rogers--Ramanujan type) partitions where the gap between each consecutive part is $\geq 2$.  Then,  \begin{equation}\label{Ramanujan_F_LHS}\sum_{n\geq 0} \frac{q^{n^2} x^n}{(q;q)_n} = \sum_{\pi\in\R\R} x^{\nu(\pi)} q^{|\pi|}. \end{equation} 

To give more context to an unfamiliar reader, we start with the simple observation \[n^2 = (2n-1) + (2n-3) + \dots + 3+1.\] For a fixed positive $n$, we interpret $q^{n^2}x^n$ as it is representing the partition \[\pi_n^* = (2n-1,  2n-3,  \dots ,3,1)\in \R\R,\] and $x^n$ keeping track of the number of parts. Observe that $\pi_n^*$ is the partition into exactly $n$-parts with the smallest possible size in $\R\R$. We will call such partitions with the smallest possible size for a fixed length \textit{minimal configuration}. Note that $1/(q;q)_n $ is the generating function for the partitions $\pi'$ into $\leq n $ parts. The partitions $\pi_n^*$ and a $\pi'$ can be put together bijectively by inserting columns of the Ferrers diagram of $\pi'$ in $\pi^*_n$'s Ferrers diagram making sure that each column of the Ferrers diagram to the right of the inserted column has a smaller height. This is until that column reaches a same or higher height column in the Ferrers diagram of $\pi_n^*$. This operation is exemplified in Table~\ref{Column_insertion}. 

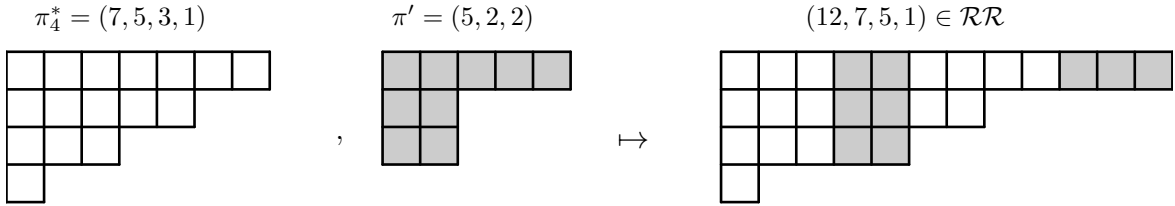
\begin{table}[htp]\caption{Example of the column insertion operation.}\label{Column_insertion}
\begin{center}
\begin{tikzpicture}[line cap=round,line join=round,x=0.5cm,y=0.5cm]
\clip(1,0.5) rectangle (32,6.5);
\fill[line width=0.pt,fill=black,fill opacity=0.2] (11.,5.) -- (16.,5.) -- (16.,4.) -- (13.,4.) -- (13.,2.) -- (11.,2.) -- cycle;
\fill[line width=0.pt,fill=black,fill opacity=0.2] (25.,2.) -- (25.,5.) -- (23.,5.) -- (23.,2.) -- cycle;
\fill[line width=0.pt,fill=black,fill opacity=0.2] (32.,5.) -- (29.,5.) -- (29.,4.) -- (32.,4.) -- cycle;
\draw [line width=1.pt] (20.,5.)-- (32.,5.);
\draw [line width=1.pt] (32.,5.)-- (32.,4.);
\draw [line width=1.pt] (32.,4.)-- (20.,4.);
\draw [line width=1.pt] (20.,3.)-- (27.,3.);
\draw [line width=1.pt] (27.,3.)-- (27.,5.);
\draw [line width=1.pt] (20.,5.)-- (20.,1.);
\draw [line width=1.pt] (20.,1.)-- (21.,1.);
\draw [line width=1.pt] (21.,1.)-- (21.,5.);
\draw [line width=1.pt] (20.,2.)-- (25.,2.);
\draw [line width=1.pt] (25.,2.)-- (25.,5.);
\draw [line width=1.pt] (26.,5.)-- (26.,3.);
\draw [line width=1.pt] (28.,5.)-- (28.,4.);
\draw [line width=1.pt] (29.,5.)-- (29.,4.);
\draw [line width=1.pt] (30.,5.)-- (30.,4.);
\draw [line width=1.pt] (31.,4.)-- (31.,5.);
\draw [line width=1.pt] (24.,5.)-- (24.,2.);
\draw [line width=1.pt] (23.,2.)-- (23.,5.);
\draw [line width=1.pt] (22.,5.)-- (22.,2.);
\draw [line width=1.pt] (11.,5.)-- (16.,5.);
\draw [line width=1.pt] (16.,5.)-- (16.,4.);
\draw [line width=1.pt] (16.,4.)-- (11.,4.);
\draw [line width=1.pt] (11.,5.)-- (11.,2.);
\draw [line width=1.pt] (11.,2.)-- (13.,2.);
\draw [line width=1.pt] (13.,2.)-- (13.,5.);
\draw [line width=1.pt] (15.,5.)-- (15.,4.);
\draw [line width=1.pt] (14.,5.)-- (14.,4.);
\draw [line width=1.pt] (12.,5.)-- (12.,2.);
\draw [line width=1.pt] (11.,3.)-- (13.,3.);
\draw [line width=1.pt] (1.,5.)-- (8.,5.);
\draw [line width=1.pt] (8.,5.)-- (8.,4.);
\draw [line width=1.pt] (8.,4.)-- (1.,4.);
\draw [line width=1.pt] (1.,3.)-- (6.,3.);
\draw [line width=1.pt] (1.,5.)-- (1.,1.);
\draw [line width=1.pt] (1.,1.)-- (2.,1.);
\draw [line width=1.pt] (2.,1.)-- (2.,5.);
\draw [line width=1.pt] (1.,2.)-- (4.,2.);
\draw [line width=1.pt] (4.,2.)-- (4.,5.);
\draw [line width=1.pt] (6.,3.)-- (6.,5.);
\draw [line width=1.pt] (7.,5.)-- (7.,4.);
\draw [line width=1.pt] (5.,5.)-- (5.,3.);
\draw [line width=1.pt] (3.,2.)-- (3.,5.);
\draw (1.5,6.5) node[anchor=north west] {$\pi_4^* = (7,5,3,1)$};
\draw (11,6.5) node[anchor=north west] {$\pi' = (5,2,2)$};
\draw (22,6.5) node[anchor=north west] {$(12,7,5,1)\in\R\R$};
\draw (17,3) node[anchor=north west] {\Large $\mapsto$};
\draw (9.5,3) node[anchor=north west] {\Large $,$};
\end{tikzpicture}
\end{center}
\end{table}

Also one should observe that the inverse operation of taking columns out of any arbitrary partition, that belongs to $\R\R$, into exactly $n$ parts to make a pair of partitions where one partition is the minimal configuration $\pi_n^*$ and another partition (may be the empty partition) into $\leq n$ parts is straightforward. Therefore, put together \[\frac{q^{n^2}x^n}{(q;q)_n} =\sum_{\substack{\pi\in\R\R\\\nu(\pi)=n}} x^{\nu(\pi)} q^{|\pi|}. \] Summing over all $n\geq 0$ gives us the interpretation represented in \eqref{Ramanujan_F_LHS}. 
 
The right-hand side of \eqref{Ramanujan_F} can also be interpreted similarly. For any fixed non-negative $n$, the factor \[\frac{q^{n^2}x^n}{(q^2;q^2)_n}\] of the summand is the generating function for the number of partitions into $n$ distinct odd parts, and \[(-xq^{2n+2};q^2)_\infty\] is the generating function for the number of partitions into distinct even parts $\geq 2n+2$, where $x$ is counting the number or parts in both generating functions. Any two partitions counted by these two different factors can be put together bijectively (and easily, due to the parity difference) by just joining the partitions and reordering the sequence in non-increasing order. This leads to the following interpretation  \begin{equation}\label{Ramanujan_F_RHS} \sum_{n\geq 0}\frac{q^{n^2}x^n}{(q^2;q^2)_n} (-xq^{2n+2};q^2)_\infty = \sum_{\substack{\pi\in\D\\ s_e(\pi) > 2 \nu_o(\pi)}} x^{\nu(\pi)} q^{|\pi|}, \end{equation} where $s_e(\pi)$ is the smallest even part in $\pi$ (if exists, otherwise we can define it to be $\infty$), and $\nu_o(\pi)$ is the number of odd parts in $\pi$.

	Comparing \eqref{Ramanujan_F_LHS} and \eqref{Ramanujan_F_RHS} proves Theorem~\ref{Thm_Ramanujan_F_combinatorial}.

Now we move on to the proof of \eqref{Ramanujan_F} by showing that the right-hand side series can be interpreted as the generating function counting the same partitions as the left-hand side. 	Here we would like to study the generating functions more carefully to explicitly show that the right-hand side \eqref{Ramanujan_F} can also be interpreted as generating Rogers--Ramanujan type partitions starting from unique minimal configurations.

 With this motivation, we would like to rewrite the term $(-xq^{2n+2};q^2)_\infty$ as a series first. For a fixed non-negative $n$, a partition $\pi$ counted by $(-xq^{2n+2};q^2)_\infty$ has even parts $\geq 2n+2$. If $\pi$ has $k$-parts, one can extract (subtract) $2n$ from each and every part of $\pi$ and still have a partition $\pi'$ into $k$ distinct even parts. Moreover, from $\pi'$ one can group out $2{k+1\choose 2} = k(k+1)$ more by extracting 2 from the smallest part, 4 from the second smallest, 6 from the third, ..., $2k$ from the largest part of $\pi$. This is possible since the parts are all distinct and even. This leaves a partition $\pi''$, which is a partition into $\leq k$ even parts (possibly empty). The extractions here are all bijective and can easily be reversed. Therefore, we can put together the generating functions for such partitions simply by multiplying them. The generating function for all such $\pi''$ is $1/(q^2;q^2)_k$, and the extractions can be expressed as the $q$-factor $q^{2nk + k(k+1)}$. Hence, \[\frac{q^{2nk + k(k+1)}}{(q^2;q^2)_k}\] is the generating function for the partitions into exactly $k$ distinct even parts all greater than or equal to $2n+2$. Summing over all possible $k$, and adding the variable $x^k$ to count the number of parts, it is clear now that \begin{equation}\label{dist_even_learger_than_n}(-xq^{2n+2};q^2)_\infty = \sum_{k\geq 0} \frac{x^k q^{2nk + k(k+1)}}{(q^2;q^2)_k}.\end{equation} We write \eqref{dist_even_learger_than_n} in on the right-hand side of \eqref{Ramanujan_F}, and rearrange the sums, and introduce a $q$-binomial factor in the spirit of \eqref{Main_Connection}: 
\begin{align}
\nonumber\sum_{n\geq 0}\frac{q^{n^2}x^n}{(q^2;q^2)_n} (-xq^{2n+2};q^2)_\infty &= \sum_{n\geq 0}\sum_{k\geq 0}\frac{q^{n^2+2nk+k(k+1)}x^{n+k}}{(q^2;q^2)_n(q^2;q^2)_k} \\                                                                                                                                                                                                                                                                                                                                                                                                                                                                                                                                                                                                                                                                                                                                                                                                                                                                                                                                                                                                                                                                                                                                                                                                                                                                                                                                                                                                                                                                                                                                                                                                                                                                                                                                                                                                                     
\label{n+k_representation}&= \sum_{n\geq 0}\sum_{k\geq 0}\frac{q^{n^2+2nk+k(k+1)}x^{n+k}}{(q^2;q^2)_{n+k}}{n+k\brack k}_{q^2}                                                                                                                                                                                                                                                                                                                                                                                                                                                                                                                                                                                                                                                                                                                                                                                                                                                                                                                                                                                                                                                                                                                                                                                                                                                                                                                                                                                                                                                                                                                                                                                                                       \end{align}

Taking a peek at our goal, we would like to show that the factor \[q^{n^2+2nk+k(k+1)}{n+k\brack k}_{q^2} \] is the generating function for the number of all partitions that satisfy the gap condition of $\R\R$ into $n$-odd and $k$-even parts where one can not make any column extractions without violating the $\R\R$ gap conditions. We will call such partitions {\it minimal partitions}. Hence, the extra $1/(q^2;q^2)_{n+k}$ factor introduces the extra columns to insert (as in Table~\ref{Column_insertion}) without changing parities of parts. 
The definition of minimality here is analogous to the previous case; a minimal partition from $\R\R$ is a partition of the smallest size that satisfies the gap condition of $\R\R$.
This is equivalent to saying that one can not do any column extractions from a minimal partition from $\R\R$ without violating the gap conditions. In general, for any given partition $\pi$ in $\R\R$ into $m$ parts one can bijectively extract columns of 2-thickness from as many times they can without violating the gap conditions of $\R\R$ and create a pair of partitions $(\pi_1,\pi_2)$, where $\pi_1$ is minimal and $\pi_2$ is made up of $\leq m$ even parts. The generating function for all such $\pi_2$ is \begin{equation}\label{Even_column_insertions_GF}
\frac{1}{(q^2;q^2)_M}.
\end{equation}

We still need to study the minimal partitions from $\R\R$ more and make some observations. For visualization, in Table~\ref{Most_minimal_n_k}, we would like to start with the interpretation of $q^{n^2+2nk+k(k+1)}$ as a 2-modular Ferrers diagram \cite[Ex 5., p. 13]{A}. (The terms in the exponent can be written more elegantly but we will leave these three terms separate to keep the connection with extractions outlined in \eqref{dist_even_learger_than_n}.)

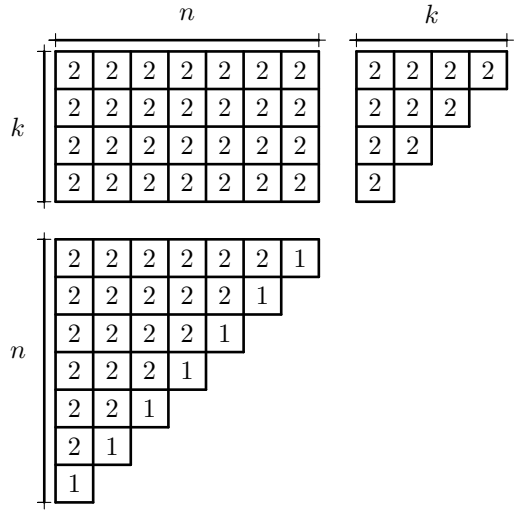
\begin{table}[thp]\caption{The visualization of the factor $q^{n^2+2nk+k(k+1)}$ as a partition in $\R\R$ using $(n,k) = (7,4)$.}\label{Most_minimal_n_k}
\begin{center}
\begin{tikzpicture}[line cap=round,line join=round,>=triangle 45,x=0.5cm,y=0.5cm]
\clip(0.5,0.5) rectangle (14.5,15.);
\draw [line width=1.pt] (2.,13.)-- (9.,13.);
\draw [line width=1.pt] (9.,13.)-- (9.,9.);
\draw [line width=1.pt] (9.,9.)-- (2.,9.);
\draw [line width=1.pt] (2.,9.)-- (2.,13.);
\draw [line width=1.pt] (2.,12.)-- (9.,12.);
\draw [line width=1.pt] (9.,11.)-- (2.,11.);
\draw [line width=1.pt] (2.,10.)-- (9.,10.);
\draw [line width=1.pt] (10.,9.)-- (11.,9.);
\draw [line width=1.pt] (11.,9.)-- (11.,13.);
\draw [line width=1.pt] (10.,13.)-- (14.,13.);
\draw [line width=1.pt] (14.,13.)-- (14.,12.);
\draw [line width=1.pt] (14.,12.)-- (10.,12.);
\draw [line width=1.pt] (10.,13.)-- (10.,9.);
\draw [line width=1.pt] (10.,10.)-- (12.,10.);
\draw [line width=1.pt] (10.,11.)-- (12.,11.);
\draw [line width=1.pt] (12.,13.)-- (12.,10.);
\draw [line width=1.pt] (12.,11.)-- (13.,11.);
\draw [line width=1.pt] (13.,11.)-- (13.,13.);
\draw [line width=1.pt] (8.,13.)-- (8.,9.);
\draw [line width=1.pt] (7.,9.)-- (7.,13.);
\draw [line width=1.pt] (6.,13.)-- (6.,9.);
\draw [line width=1.pt] (5.,9.)-- (5.,13.);
\draw [line width=1.pt] (4.,13.)-- (4.,9.);
\draw [line width=1.pt] (3.,9.)-- (3.,13.);
\draw [line width=1.pt] (9.,8.)-- (2.,8.);
\draw [line width=1.pt] (2.,8.)-- (2.,1.);
\draw [line width=1.pt] (2.,1.)-- (3.,1.);
\draw [line width=1.pt] (3.,1.)-- (3.,8.);
\draw [line width=1.pt] (9.,7.)-- (9.,8.);
\draw [line width=1.pt] (8.,6.)-- (8.,8.);
\draw [line width=1.pt] (7.,5.)-- (7.,8.);
\draw [line width=1.pt] (6.,4.)-- (6.,8.);
\draw [line width=1.pt] (5.,3.)-- (5.,8.);
\draw [line width=1.pt] (4.,2.)-- (4.,8.);
\draw [line width=1.pt] (9.,7.)-- (2.,7.);
\draw [line width=1.pt] (8.,6.)-- (2.,6.);
\draw [line width=1.pt] (7.,5.)-- (2.,5.);
\draw [line width=1.pt] (6.,4.)-- (2.,4.);
\draw [line width=1.pt] (5.,3.)-- (2.,3.);
\draw [line width=1.pt] (4.,2.)-- (2.,2.);
\draw (2.5,12.5) node[anchor=center] {2};
\draw (2.5,2.5) node[anchor=center] {2};
\draw (10.5,9.5) node[anchor=center] {2};
\draw (2.5,1.5) node[anchor=center] {1};
\draw [line width=1.pt] (2.,13.3)-- (9.,13.3);
\draw [line width=1.pt] (10.,13.3)-- (14.,13.3);
\draw [line width=1.pt] (1.7,13.)-- (1.7,9.);
\draw [line width=1.pt] (1.7,8.)-- (1.7,1.);
\draw (1,5) node[anchor=center] {$n$};
\draw (1,11) node[anchor=center] {$k$};
\draw (5.5,14) node[anchor=center] {$n$};
\draw (12,14) node[anchor=center] {$k$};
\draw (3.5,2.5) node[anchor=center] {1};
\draw (4.5,3.5) node[anchor=center] {1};
\draw (5.5,4.5) node[anchor=center] {1};
\draw (6.5,5.5) node[anchor=center] {1};
\draw (7.5,6.5) node[anchor=center] {1};
\draw (8.5,7.5) node[anchor=center] {1};
\draw (3.5,12.5) node[anchor=center] {2};
\draw (4.5,12.5) node[anchor=center] {2};
\draw (5.5,12.5) node[anchor=center] {2};
\draw (6.5,12.5) node[anchor=center] {2};
\draw (7.5,12.5) node[anchor=center] {2};
\draw (8.5,12.5) node[anchor=center] {2};
\draw (2.5,11.5) node[anchor=center] {2};
\draw (3.5,11.5) node[anchor=center] {2};
\draw (4.5,11.5) node[anchor=center] {2};
\draw (5.5,11.5) node[anchor=center] {2};
\draw (6.5,11.5) node[anchor=center] {2};
\draw (7.5,11.5) node[anchor=center] {2};
\draw (8.5,11.5) node[anchor=center] {2};
\draw (2.5,10.5) node[anchor=center] {2};
\draw (3.5,10.5) node[anchor=center] {2};
\draw (4.5,10.5) node[anchor=center] {2};
\draw (5.5,10.5) node[anchor=center] {2};
\draw (6.5,10.5) node[anchor=center] {2};
\draw (7.5,10.5) node[anchor=center] {2};
\draw (8.5,10.5) node[anchor=center] {2};
\draw (2.5,9.5) node[anchor=center] {2};
\draw (3.5,9.5) node[anchor=center] {2};
\draw (4.5,9.5) node[anchor=center] {2};
\draw (5.5,9.5) node[anchor=center] {2};
\draw (6.5,9.5) node[anchor=center] {2};
\draw (7.5,9.5) node[anchor=center] {2};
\draw (8.5,9.5) node[anchor=center] {2};
\draw (2.5,7.5) node[anchor=center] {2};
\draw (3.5,7.5) node[anchor=center] {2};
\draw (4.5,7.5) node[anchor=center] {2};
\draw (5.5,7.5) node[anchor=center] {2};
\draw (6.5,7.5) node[anchor=center] {2};
\draw (7.5,7.5) node[anchor=center] {2};
\draw (2.5,6.5) node[anchor=center] {2};
\draw (3.5,6.5) node[anchor=center] {2};
\draw (4.5,6.5) node[anchor=center] {2};
\draw (5.5,6.5) node[anchor=center] {2};
\draw (6.5,6.5) node[anchor=center] {2};
\draw (10.5,12.5) node[anchor=center] {2};
\draw (11.5,12.5) node[anchor=center] {2};
\draw (10.5,11.5) node[anchor=center] {2};
\draw (11.5,11.5) node[anchor=center] {2};
\draw (2.5,5.5) node[anchor=center] {2};
\draw (3.5,5.5) node[anchor=center] {2};
\draw (2.5,4.5) node[anchor=center] {2};
\draw (3.5,4.5) node[anchor=center] {2};
\draw (12.5,12.5) node[anchor=center] {2};
\draw (13.5,12.5) node[anchor=center] {2};
\draw (12.5,11.5) node[anchor=center] {2};
\draw (4.5,5.5) node[anchor=center] {2};
\draw (5.5,5.5) node[anchor=center] {2};
\draw (4.5,4.5) node[anchor=center] {2};
\draw (10.5,10.5) node[anchor=center] {2};
\draw (11.5,10.5) node[anchor=center] {2};
\draw (2.5,3.5) node[anchor=center] {2};
\draw (3.5,3.5) node[anchor=center] {2};
\begin{scriptsize}
\draw [color=black] (2.,13.3)-- ++(-2.5pt,0 pt) -- ++(5.0pt,0 pt) ++(-2.5pt,-2.5pt) -- ++(0 pt,5.0pt);
\draw [color=black] (9.,13.3)-- ++(-2.5pt,0 pt) -- ++(5.0pt,0 pt) ++(-2.5pt,-2.5pt) -- ++(0 pt,5.0pt);
\draw [color=black] (10.,13.3)-- ++(-2.5pt,0 pt) -- ++(5.0pt,0 pt) ++(-2.5pt,-2.5pt) -- ++(0 pt,5.0pt);
\draw [color=black] (14.,13.3)-- ++(-2.5pt,0 pt) -- ++(5.0pt,0 pt) ++(-2.5pt,-2.5pt) -- ++(0 pt,5.0pt);
\draw [color=black] (1.7,13.)-- ++(-2.5pt,0 pt) -- ++(5.0pt,0 pt) ++(-2.5pt,-2.5pt) -- ++(0 pt,5.0pt);
\draw [color=black] (1.7,9.)-- ++(-2.5pt,0 pt) -- ++(5.0pt,0 pt) ++(-2.5pt,-2.5pt) -- ++(0 pt,5.0pt);
\draw [color=black] (1.7,8.)-- ++(-2.5pt,0 pt) -- ++(5.0pt,0 pt) ++(-2.5pt,-2.5pt) -- ++(0 pt,5.0pt);
\draw [color=black] (1.7,1.)-- ++(-2.5pt,0 pt) -- ++(5.0pt,0 pt) ++(-2.5pt,-2.5pt) -- ++(0 pt,5.0pt);
\end{scriptsize}
\end{tikzpicture}
\end{center}
\end{table}

The partition encoded by $q^{n^2+2nk+k(k+1)}$ ---with no excess columns that can be extracted--- is a minimal partition in the set $\R\R$ into $n$-odd and $k$-even parts. We want to point out that there can be multiple\footnote{The exact amount is ${n+k\choose n}.$} minimal partitions for such restrictions depending on the location of the odd parts. To demonstrate we list all 6 minimal partitions of $\R\R$ into 2-odd and 2-even parts as 2-modular Ferrers diagrams in Table~\ref{Minimal_2_2}.

\begin{table}[thp]\caption{Minimal partitions of $\R\R$ into 2-odd and 2-even parts as 2-modular Ferrers diagrams.}\label{Minimal_2_2}
\begin{center}
\begin{tikzpicture}[line cap=round,line join=round,x=0.5cm,y=0.5cm]
\clip(1.5,0.5) rectangle (22.5,11.5);
\draw [line width=1.pt] (2.,11.)-- (6.,11.);
\draw [line width=1.pt] (6.,11.)-- (6.,10.);
\draw [line width=1.pt] (6.,10.)-- (2.,10.);
\draw [line width=1.pt] (2.,11.)-- (2.,7.);
\draw [line width=1.pt] (2.,7.)-- (3.,7.);
\draw [line width=1.pt] (3.,7.)-- (3.,11.);
\draw [line width=1.pt] (2.,9.)-- (5.,9.);
\draw [line width=1.pt] (5.,9.)-- (5.,11.);
\draw [line width=1.pt] (4.,11.)-- (4.,8.);
\draw [line width=1.pt] (4.,8.)-- (2.,8.);
\draw [line width=1.pt] (2.,5.)-- (7.,5.);
\draw [line width=1.pt] (7.,5.)-- (7.,4.);
\draw [line width=1.pt] (7.,4.)-- (2.,4.);
\draw [line width=1.pt] (2.,5.)-- (2.,1.);
\draw [line width=1.pt] (2.,1.)-- (3.,1.);
\draw [line width=1.pt] (3.,1.)-- (3.,5.);
\draw [line width=1.pt] (2.,3.)-- (6.,3.);
\draw [line width=1.pt] (6.,3.)-- (6.,5.);
\draw [line width=1.pt] (5.,5.)-- (5.,3.);
\draw [line width=1.pt] (4.,5.)-- (4.,2.);
\draw [line width=1.pt] (4.,2.)-- (2.,2.);
\draw [line width=1.pt] (9.,1.)-- (10.,1.);
\draw [line width=1.pt] (9.,1.)-- (9.,5.);
\draw [line width=1.pt] (9.,5.)-- (14.,5.);
\draw [line width=1.pt] (14.,5.)-- (14.,4.);
\draw [line width=1.pt] (14.,4.)-- (9.,4.);
\draw [line width=1.pt] (9.,3.)-- (13.,3.);
\draw [line width=1.pt] (11.,2.)-- (9.,2.);
\draw [line width=1.pt] (10.,1.)-- (10.,5.);
\draw [line width=1.pt] (11.,2.)-- (11.,5.);
\draw [line width=1.pt] (13.,3.)-- (13.,5.);
\draw [line width=1.pt] (12.,5.)-- (12.,3.);
\draw [line width=1.pt] (17.,1.)-- (16.,1.);
\draw [line width=1.pt] (16.,1.)-- (16.,5.);
\draw [line width=1.pt] (16.,5.)-- (22.,5.);
\draw [line width=1.pt] (22.,5.)-- (22.,4.);
\draw [line width=1.pt] (22.,4.)-- (16.,4.);
\draw [line width=1.pt] (16.,3.)-- (20.,3.);
\draw [line width=1.pt] (21.,5.)-- (21.,4.);
\draw [line width=1.pt] (20.,5.)-- (20.,3.);
\draw [line width=1.pt] (19.,5.)-- (19.,2.);
\draw [line width=1.pt] (19.,2.)-- (16.,2.);
\draw [line width=1.pt] (17.,1.)-- (17.,5.);
\draw [line width=1.pt] (18.,5.)-- (18.,2.);
\draw [line width=1.pt] (17.,7.)-- (16.,7.);
\draw [line width=1.pt] (16.,7.)-- (16.,11.);
\draw [line width=1.pt] (16.,11.)-- (21.,11.);
\draw [line width=1.pt] (21.,11.)-- (21.,10.);
\draw [line width=1.pt] (21.,10.)-- (16.,10.);
\draw [line width=1.pt] (20.,9.)-- (16.,9.);
\draw [line width=1.pt] (18.,8.)-- (16.,8.);
\draw [line width=1.pt] (17.,7.)-- (17.,11.);
\draw [line width=1.pt] (18.,11.)-- (18.,8.);
\draw [line width=1.pt] (19.,11.)-- (19.,9.);
\draw [line width=1.pt] (20.,9.)-- (20.,11.);
\draw [line width=1.pt] (9.,11.)-- (14.,11.);
\draw [line width=1.pt] (14.,11.)-- (14.,10.);
\draw [line width=1.pt] (14.,10.)-- (9.,10.);
\draw [line width=1.pt] (9.,11.)-- (9.,7.);
\draw [line width=1.pt] (9.,7.)-- (10.,7.);
\draw [line width=1.pt] (10.,7.)-- (10.,11.);
\draw [line width=1.pt] (9.,9.)-- (12.,9.);
\draw [line width=1.pt] (12.,9.)-- (12.,11.);
\draw [line width=1.pt] (11.,11.)-- (11.,8.);
\draw [line width=1.pt] (11.,8.)-- (9.,8.);
\draw [line width=1.pt] (13.,11.)-- (13.,10.);
\draw (2.5,7.5) node[anchor=center] {$1$};
\draw (3.5,8.5) node[anchor=center] {$1$};
\draw (13.5,10.5) node[anchor=center] {$1$};
\draw (9.5,7.5) node[anchor=center] {$1$};
\draw (19.5,9.5) node[anchor=center] {$1$};
\draw (20.5,10.5) node[anchor=center] {$1$};
\draw (5.5,3.5) node[anchor=center] {$1$};
\draw (2.5,1.5) node[anchor=center] {$1$};
\draw (12.5,3.5) node[anchor=center] {$1$};
\draw (18.5,2.5) node[anchor=center] {$1$};
\draw (21.5,4.5) node[anchor=center] {$1$};
\draw [line width=1.pt] (12.,2.)-- (12.,3.);
\draw [line width=1.pt] (12.,2.)-- (11.,2.);
\draw (11.5,2.5) node[anchor=center] {$1$};
\draw (2.5,10.5) node[anchor=center]{$2$};
\draw (3.5,10.5) node[anchor=center]{$2$};
\draw (4.5,10.5) node[anchor=center]{$2$};
\draw (5.5,10.5) node[anchor=center]{$2$};
\draw (4.5,9.5) node[anchor=center]{$2$};
\draw (3.5,9.5) node[anchor=center]{$2$};
\draw (2.5,9.5) node[anchor=center]{$2$};
\draw (2.5,8.5) node[anchor=center]{$2$};
\draw (9.5,10.5) node[anchor=center]{$2$};
\draw (10.5,10.5) node[anchor=center]{$2$};
\draw (11.5,10.5) node[anchor=center]{$2$};
\draw (12.5,10.5) node[anchor=center]{$2$};
\draw (11.5,9.5) node[anchor=center]{$2$};
\draw (10.5,9.5) node[anchor=center]{$2$};
\draw (9.5,9.5) node[anchor=center]{$2$};
\draw (9.5,8.5) node[anchor=center]{$2$};
\draw (16.5,10.5) node[anchor=center]{$2$};
\draw (17.5,10.5) node[anchor=center]{$2$};
\draw (18.5,10.5) node[anchor=center]{$2$};
\draw (19.5,10.5) node[anchor=center]{$2$};
\draw (18.5,9.5) node[anchor=center]{$2$};
\draw (17.5,9.5) node[anchor=center]{$2$};
\draw (16.5,9.5) node[anchor=center]{$2$};
\draw (16.5,8.5) node[anchor=center]{$2$};
\draw (2.5,4.5) node[anchor=center]{$2$};
\draw (3.5,4.5) node[anchor=center]{$2$};
\draw (4.5,4.5) node[anchor=center]{$2$};
\draw (5.5,4.5) node[anchor=center]{$2$};
\draw (4.5,3.5) node[anchor=center]{$2$};
\draw (3.5,3.5) node[anchor=center]{$2$};
\draw (2.5,3.5) node[anchor=center]{$2$};
\draw (2.5,2.5) node[anchor=center]{$2$};
\draw (9.5,4.5) node[anchor=center]{$2$};
\draw (10.5,4.5) node[anchor=center]{$2$};
\draw (11.5,4.5) node[anchor=center]{$2$};
\draw (12.5,4.5) node[anchor=center]{$2$};
\draw (11.5,3.5) node[anchor=center]{$2$};
\draw (10.5,3.5) node[anchor=center]{$2$};
\draw (9.5,3.5) node[anchor=center]{$2$};
\draw (9.5,2.5) node[anchor=center]{$2$};
\draw (16.5,4.5) node[anchor=center]{$2$};
\draw (17.5,4.5) node[anchor=center]{$2$};
\draw (18.5,4.5) node[anchor=center]{$2$};
\draw (19.5,4.5) node[anchor=center]{$2$};
\draw (18.5,3.5) node[anchor=center]{$2$};
\draw (17.5,3.5) node[anchor=center]{$2$};
\draw (16.5,3.5) node[anchor=center]{$2$};
\draw (16.5,2.5) node[anchor=center]{$2$};
\draw (10.5,8.5) node[anchor=center]{$2$};
\draw (16.5,7.5) node[anchor=center]{$2$};
\draw (17.5,8.5) node[anchor=center]{$2$};
\draw (19.5,3.5) node[anchor=center]{$2$};
\draw (20.5,4.5) node[anchor=center]{$2$};
\draw (17.5,2.5) node[anchor=center]{$2$};
\draw (16.5,1.5) node[anchor=center]{$2$};
\draw (9.5,1.5) node[anchor=center]{$2$};
\draw (10.5,2.5) node[anchor=center]{$2$};
\draw (13.5,4.5) node[anchor=center]{$2$};
\draw (3.5,2.5) node[anchor=center]{$2$};
\draw (6.5,4.5) node[anchor=center]{$2$};
\end{tikzpicture}
\end{center}
\end{table}

Now, to simplify the notation, we would like to pair minimal partitions with vectors. There is a bijection between the minimal partitions in $\R\R$ with $n$-odd and $k$-even parts and the $1\times(n+k)$-vectors that has $k$ $1$s and $n$ $2$s as their entries. From the minimal partitions to the vectors is easy to express; Starting from the smallest part (bottom row of the Ferrers diagram) we put the end-of-row box value of the 2-modular Ferrers diagram in the entries of the vector. This is exemplified in Table~\ref{Vectors_2_2}. 

\begin{table}[thp]\caption{Respective vector representations of the minimal partitions in $\R\R$ into 2-odd and 2-even parts presented in Table~\ref{Minimal_2_2}}\label{Vectors_2_2}
\[\begin{array}{c}
(1,1,2,2),\quad (1,2,2,1),\quad (2,2,1,1), \\
(1,2,1,2),\quad (2,1,1,2),\quad (2,1,2,1).
\end{array}\]
\end{table}

The inverse construction, from a vector to a 2-modular Ferrers diagram/partition, needs us to make some observations. First observation is that the partition encoded by $q^{n^2+2nk+k(k+1)}$ is \[(\underbrace{1,1,\dots,1}_{n-\text{many}},\underbrace{2,2,\dots,2}_{k-\text{many}}).\]
One can construct a 2-modular Ferrers diagram from a given vector of 1s and 2s. Let $M$ be the length of such a vector $v$, one starts by creating a triangular Ferrers diagram $\pi^*$ of height $M$, fills the ends of each column with the vector's entries and the rest of the boxes (that lie in the interior of the 2-modular Ferrers diagram) with 2s. This Ferrers diagram $\pi^*$ represents a partition, but it might be violating the difference condition of $\R\R$ between the parts. For the gap conditions then one inserts an intermediate column of boxes filled with 2s in the respective places where the vector of $\pi^*$ has the pair $(\dots, 2,1, \dots );$ a $2$ followed by a $1$. The height of this intermediate column is determined by the index of the related $1$ in the vector; if the index of the 1 is in the $i$-th position then the related column height is $M-i+1$. Observe that the extra columns that needs to be inserted all have the height $\leq M-1$, as we never need to insert a column before the first column, which has the height $M$. Also, observe that each extra column's height differs by each other by at least 2.

Now that all these observations are laid out, we would like to write the generating function of all such minimal partitions in $\R\R$. Let $\pi$ be a minimal partition in $\R\R$  into $n$-odd and $k$-even parts, and let $v$ be its vector representation. We would like to create a pair of partitions $(\pi_1,\pi_2)$ from $\pi$. We extract any intermediate columns that $\pi$ has and put those in $\pi_2$ and we denote the leftover partition as $\pi_1$. Unsurprisingly, $\pi_1$ has size $q^{n^2+2nk+k(k+1)}$.

Our remaining task is to show that the generating function for all such $\pi_2$ is \[{n+k \brack k}_{q^2}.\] Equivalently, we would like to show that $\pi_2$ and a given minimal configuration vector $v$ is in bijection with 2-modular Ferrers diagrams that fit in a $n\times k$ box. We will use the ideas behind Sylvester's bijection \cite{BU2} for this remaining task and explain, algorithmically, how one can pair partitions $\pi_2$ using the vector $v$, and partitions that fit in a $n\times k$ box, where all the boxes are filled with 2s. 

\begin{enumerate}[i.]
\item One takes the shortest column of the 2-modular Ferrers diagram of $\pi_2$. The height $h$ of this column is related with a 1 in the vector $v$ at the index $(n+k) - h + 1$. There are $n_h$ 1s and $k_h$ 2s from this index to the end of the vector.
\item One takes the $h$ boxes in the selected column and puts $n_h$ of them in vertical orientation and then puts the remaining $k_h$ of these boxes in a horizontal orientation right next to the top box in the $n_h$-height column. This creates a hook made up of a vertical \textit{leg} (of size $n_h$) and a horizontal \textit{arm} (of size $k_h +1$).
\item Then one repeats the above construction as in Sylvester's bijection and puts the created hook on top of the previous hook pile to create a partition $\hat{\pi}$.
\end{enumerate}

%This operation is exemplified in Table~\eqref{Table_long_example}. 

\begin{table}

\end{table}

One needs to convince themselves that $\hat{\pi}$ is indeed a partition. We want to point out that, for us to get a partition from this construction, both the leg and arm of a hook must be at least one longer than of the previous hook's leg and arm. A new hook gets created if a 2 is followed by a 1 in the vector $v$ by using that said 1 and the contents of $v$ untill the end. Usage of the hook-creating 1 shows us that there is at least one extra 1 in the system and the leg of that hook will be longer, and the 2 that created the previous hook is now in consideration and this similarly shows that the arm is longer. This is hand-in-hand with the column heights of $\pi_2$ being at least 2 different from each other. 

For fixed $n$ and $k$, the last two points to be addressed are the outcome partition's size (and if it fits in a $n\times k$-box), and the bijectivity of this map (and whether we can find a $v$ and a $\pi_2$ given a partition in a $n\times k$-box). The largest (last) hook always has $\geq n$ leg length and $\leq k$ arm length, since even at the extreme case (where a 2 is in the first index of the vector $v$ followed by a 1) one can take all 1s in consideration but all-but-one of the 2s. Hence, the created partition fits in a $n\times k$ box. Given a partition $\pi$ that fits in a $n\times k$ box, unfolding the hooks to make the partition $\pi_2$ is trivial. Finally, to create the vector $v$ one starts from $\pi$ we need to look at the hooks. If $\pi$ is the empty partition we are in the base case and $v$ is made up of $n$ 1s followed by $k$ 2s. If $\pi$ is not empty, we look at the innermost hook and write leg-length many 1s followed by arm-length minus 1 many 2s in a vector. Then, we move on to the next hook, and prior to what we have already written down we write leg-length many 1s followed by arm-length minus 1 many 2s. We repeat this procedure till we run out of hooks. This creates a vector $v'$ of length $M=m_1+m_2$ with $m_1$ 1s and $m_2$ 2s. We create the expected representative vector $v$ by appending $(n-m_1)$-many 1s followed by $(k-m_2)$-many 2s in front of $v'$. This shows that there is a bijection between the vector, partition pair $(v,\pi_2)$, and partitions that fit in a $n\times k$ box.

This study shows that the generating function for the number of minimal partitions in $\R\R$ with $n$-odd and $k$-even parts is \[q^{n^2+2nk+k(k+1)}{n+k\brack k}_{q^2}. \]
Then thinking of column insertions, and keeping \eqref{Even_column_insertions_GF} in mind \[\frac{q^{n^2+2nk+k(k+1)}}{(q^2;q^2)_{n+k}}{n+k\brack k}_{q^2} \] is the generating function for the number of all $\R\R$ partitions that has $n$-odd and $k$-even parts. Hence, summing over all possible $n$ and $k$ and keeping track of the total number of parts as the exponent of $x$ we see that

\begin{equation}\label{RR_last}
\sum_{n\geq 0}\sum_{k\geq 0}\frac{q^{n^2+2nk+k(k+1)}x^{n+k}}{(q^2;q^2)_{n+k}}{n+k\brack k}_{q^2}                                                                                    = \sum_{\pi\in\R\R} x^{\nu(\pi)} q^{|\pi|}.
\end{equation}

Hence, comparing \eqref{RR_last} and \eqref{n+k_representation}, this shows that the left-hand side of \eqref{n+k_representation} is equal to the right-hand side of \eqref{Ramanujan_F_LHS}. This finishes the proof of the identity \eqref{Ramanujan_F} and implies Theorem~\ref{Thm2}.

\section{Outlook}

The $q$-series identities that came out of Mourad's studies seem to be receptive to combinatorial constructions. It would be interesting to see if more identities can be proven this way. Moreover, it would be of interest to see if new identities that were not discovered by orthogonal polynomial research can be found by first interpreting an object of interest and then using partition theoretic bijective modifications.

The combinatorial constructions used in these proofs resemble the authors' and University of Florida Number theory group's earlier works \cite{Alladi1, BD1, D, BU1, BU2, Uncu1}, especially \cite{Uncu2}. It should be noted that in most of these works, once the combinatorial structure behind an identity is understood, the author and Berkovich introduced extra bounds on the partitions and refined the $q$-series identities to polynomial identities that depend on one or more discrete variables. It is no different here. Theorem~\ref{Thm4} is the refinement of comparing the left-hand sides of \eqref{RR_last} and \ref{Ramanujan_F_LHS}. We study how many columns can be added to odd ($n$-many) and even ($k$-many) parts of a minimal configuration that is defined in Section~\ref{S2}, and limit these contributions to assure all the parts stay less than or equal to $L$. This identity can easily be proven by $q$-Zeilberger algorithm \cite{Z}. We omit it here. 

Such refined identities would benefit orthogonal polynomial research as well as combinatorics research. For example, knowing Theorem~\ref{Thm3}, we can experiment, conjecture the following:

\begin{theorem}\label{Thm3refined}
	\[\sum_{\substack{(\pi_1,\pi_2)\in \D_N\times\P_N\\l(\pi_2)\leq \nu(\pi_1)}} t^{\nu(\pi_1)+\nu(\pi_2)}a^{\nu(\pi_2)}b^{\nu(\pi_1)}q^{|\pi_1|+|\pi_2|} = \sum_{\substack{(\pi^*_1,\pi^*_2)\in\D_N\times\P_N\\d(\pi^*_2)\leq r(\pi^*_1)}} t^{\nu(\pi^*_1)+\nu(\pi^*_2)}a^{\nu(\pi^*_2)}b^{\nu(\pi^*_1)}q^{|\pi^*_1|+|\pi^*_2|},\] where $\D_N$ and $\P_N$ are the set of partitions from $\D$ and $\P$, respectively, with the extra condition that the parts are at most $N$.
\end{theorem}

The $q$-series equivalent of the above theorem, which we can prove easily is as follows:

\begin{theorem}\label{ThmL} For any non-negative integer $N$, we have
	\begin{equation*}\sum_{k\geq 0 }\frac{(bt)^kq^{k(k+1)/2}}{(atq;q)_k} {N\brack k}_q = \sum_{m\geq 0} (bt)^mq^{m(m+1)/2}(-btq^{m+1};q)_{N-m} \frac{(at)^mq^{m^2}}{(atq;q)_m} {N\brack m}_q.\end{equation*}
\end{theorem}

Notice that Theorem~\ref{ThmL} is a (simple) refinement of \eqref{Ismail1_simplified}, and would imply Theorem~\ref{Thm1}. We believe that Theorem~\ref{ThmL} does not appear in \cite{Mourad2019}. It would be of interest to see other proofs and implications of this theorem, whether by using bijective maps between partitions or orthogonal polynomials. In general, it is our hope to see these different techniques to bolster and push each other forward.

\section*{Acknowledgement}

The author wants to thank the anonymous referees for their careful reading and correcting the typos, expecially the ones in the formulas. Their time, effort, and rightful frustration made this paper a much better read and mistake free. The author also thanks Runqiao Li for reading the manuscript before its submission.
	
The author is grateful for the partial support of the Austrian Science Fund FWF P34501N, and UKRI EPSRC EP/T015713/1 projects.

%\section*{Conflict of Interest Statement}

%The author declares no conflict of interest.

\end{document}